\documentclass[11pt]{article}
\usepackage{listings}
\usepackage{graphicx,url}
\usepackage[utf8]{inputenc}
\usepackage{amsthm}
\usepackage{newtxtext}    
\usepackage{scrextend} 
\usepackage{soul} 
\usepackage{amsmath,amsthm,amssymb,epsfig,graphics,psfrag,latexsym,amsmath,listings,graphicx,url,tabularx,gensymb,caption,subcaption}
\usepackage{algorithm}
\usepackage[utf8]{inputenc}
\usepackage[noend]{algpseudocode}
\usepackage{multicol}
\usepackage{bbm}
\usepackage{amsrefs}
\usepackage{color}
\graphicspath{ {Figures/} }
\usepackage{geometry}
\usepackage{tikz}
\usetikzlibrary{shapes}
\usetikzlibrary{matrix}
\usetikzlibrary{positioning}
\usetikzlibrary{fit}

\tikzstyle{res}=[circle,thick,minimum size=4mm,draw=black,fill=red,inner sep=1pt]
\tikzstyle{non-res}=[circle,thick,minimum size=4mm,draw=black,inner sep=1pt]
\tikzstyle{light-res}=[circle,thick,minimum size=4mm,draw=black,fill=red!40,inner sep=1pt]
\tikzstyle{blue}=[circle,thick,minimum size=4mm,draw=black,fill=blue!20,inner sep=1pt]

\newtheorem{theorem}{Theorem}[section]
\newtheorem{corollary}{Corollary}[section]
\newtheorem{definition}{Definition}[section]
\newtheorem{lemma}{Lemma}[section]

\def\Mbf{{\mathbf{M}}}
\def\Dbf{{\mathbf{D}}}
\def\diam{{\hbox{diam}}}

\title{Truncated Metric Dimension for Finite Graphs}
\author{Richard C. Tillquist \footnote{Corresponding author}\\
        Department of Computer Science \\
        California State University, Chico \\
        \texttt{rctillquist@csuchico.edu}
        \and
        Rafael M. Frongillo \\
        Department of Computer Science \\
        University of Colorado, Boulder \\
        \texttt{raf@colorado.edu}
        \and
        Manuel E. Lladser \\
        Department of Applied Mathematics \\
        University of Colorado, Boulder \\
        \texttt{manuel.lladser@colorado.edu}}
        
\date{}

\usepackage[colorinlistoftodos,textsize=tiny]{todonotes}
\newcommand{\Comments}{1}
\newcommand{\mynote}[2]{\ifnum\Comments=1\textcolor{#1}{#2}\fi}
\newcommand{\mytodo}[2]{\ifnum\Comments=1%
  \todo[linecolor=#1!80!black,backgroundcolor=#1,bordercolor=#1!80!black]{#2}\fi}

\begin{document} 

\maketitle

\begin{abstract}
A graph $G=(V,E)$ with geodesic distance $d(\cdot,\cdot)$ is said to be resolved by a non-empty subset $R$ of its vertices when, for all vertices $u$ and $v$, if $d(u,r)=d(v,r)$ for each $r\in R$, then $u=v$. The metric dimension of $G$ is the cardinality of its smallest resolving set. In this manuscript, we present and investigate the notions of resolvability and metric dimension when the geodesic distance is truncated with a certain threshold $k$; namely, we measure distances in $G$ using the metric $d_k(u,v):=\min\{d(u,v),k+1\}$. We denote the metric dimension of $G$ with respect to $d_k$ as $\beta_k(G)$. We study the behavior of this quantity with respect to $k$ as well as the diameter of $G$. We also characterize the truncated metric dimension of paths and cycles as well as graphs with extreme metric dimension, including graphs of order $n$ such that $\beta_k(G)=n-2$ and $\beta_k(G)=n-1$. We conclude with a study of various problems related to the truncated metric dimension of trees.
\end{abstract}

\section{Introduction}

Given a simple, finite graph $G=(V,E)$, let $d(u,v)$ be the shortest path or geodesic distance between two vertices $u,v\in V$.
A non-empty subset $R\subseteq V$ is called \textit{resolving}~\cite{harary1976metric,slater1975leaves} when for all $u,v\in V$, if $d(u,r) = d(v,r)$ for all $r \in R$, then $u=v$. 
The nodes of a resolving set can be used as landmarks in a graph that allow the unique identification and representation of all vertices based on distances. In fact, $R$ resolves $G$ if and only if the vector of distances $d(u|R) := (d(u, r_1), \dots, d(u,r_{|R|}))$, for an arbitrary but fixed order of the elements of $R$, is distinct for each $u \in V$.

The \textit{metric dimension} of $G$~\cite{harary1976metric,slater1975leaves}, denoted $\beta(G)$, is the cardinality of a smallest possible resolving set of $G$.
Finding a minimal resolving set in an arbitrary graph is an NP-hard problem~\cite{gary1979computers}.
The state-of-the-art algorithm for finding non-trivial resolving sets in general graphs is based on the Information Content Heuristic (ICH)~\cite{hauptmann2012approximation} and has an approximation ratio of $1 + (1 + o(1))\cdot\ln|V|$ and a time complexity of $O(|V|^3)$.
For a short introduction to these concepts see~\cite{TilFroLla19}, and for a more comprehensive one see~\cite{tillquist2021getting}.

The traditional definition of metric dimension assumes knowledge of all pairwise distances between vertices.
This assumption allows any individual vertex $v$ of a resolving set to play a key role in distinguishing any pair of vertices, including pairs very far from $v$. In practice, however, computing pairwise distances between all pairs of vertices in a large network may be costly, and the quality of pairwise distance measurements may degrade with increasing distance.
Indeed, computing the distance matrix of a dense graph has time complexity $O(|V|^3)$~\cite{floyd1962algorithm}, whereas for a sparse graph the complexity is $O(|V||E|+|V|^2\ln|V|)$~\cite{Joh77}. 
On the other hand, various models of epidemic spread over networks assume that transmission times across edges are random but independent and identically distributed~\cite{pinto2012locating,spinelli2016observer}. As a consequence, transmission times between vertices with many intermediate edges can have high variance. In this setting, resolving sets  
using expected transmission time between pairs of vertices as the metric may not be effective in identifying the source of the infection (aka ground zero).

The above factors motivate metrics in graphs that only rely on local vertex information in the graph.
For a given integer $k\ge0$, a natural choice for such a metric is
\begin{equation}
d_k(u,v):=\min\{d(u,v),k+1\}.
\label{def:dk}
\end{equation}
Accordingly, a $k$-truncated resolving set is a non-empty set $R\subseteq V$ such that for all $u,v\in V$, if $d_k(u,r)=d_k(v,r)$ for each $r\in R$ then $u=v$. Likewise, the $k$-truncated metric dimension of $G$, which we denote $\beta_k(G)$, is defined as the cardinality of a smallest $k$-truncated resolving set.

With truncated metric dimension, elements of a resolving set are only able to distinguish vertices up to a certain distance; in particular, computation of the full distance metric is no longer necessary. In the context of epidemics, by limiting the number of relevant edges on any shortest path between elements of a resolving set and other vertices in the graph, we have better control of the uncertainty of transmission times. Indeed, since the vertices of a graph can be at most distance $k+1$ from the nearest element of a $k$-truncated resolving set, these elements tend to be spread out over the space defined by the graph.


The paper is organized as follows. After presenting formal definitions and other preliminary information in Section~\ref{sec:preliminaries}, we examine the truncated metric dimension of paths and cycles, two graph families for which traditional metric dimension is easily computed, in Section~\ref{sec:paths_cycles}. In Section~\ref{sec:extreme} we explore graph structures with 
extreme truncated metric dimension. This includes characterizations of graphs of order $n$ with truncated metric dimension $(n-2)$ and $(n-1)$. Finally, we focus on trees in Section~\ref{sec:trees}. We define conditions under which finding truncated metric dimension exactly is straightforward, present a dynamic program capable of discovering minimal resolving sets on trees when only immediate neighbors are visible, and investigate extreme constructions for graphs in this family.

\textbf{Related work.}
When $k=1$, the metric $d_1$ coincides with the so-called adjacency distance, a notion introduced in~\cite{JanOmo12} to analyze the metric dimension of lexicographic products. This metric is defined as
\[a(u,v):=\begin{cases}
0, & \text{ if $u=v$};\\
1, & \text{ if $u\sim v$};\\
2, & \text{ if $u\not\sim v$}.
\end{cases}\]
The so-called adjacency metric dimension~\cite{JanOmo12} of a graph corresponds therefore to its $1$-truncated metric dimension. 
The present paper is based on the Ph.D. thesis by the first author~\cite{Til20}, which introduced the concept of truncated metric dimension.
A review paper by all three authors~\cite{tillquist2021getting} briefly discussed truncated metric dimension.
Finally, a more recent but independent work by other authors~\cite{GenYi21} has introduced the same notion and obtained results which overlap with ours.


\section{Preliminaries}
\label{sec:preliminaries}

To formally define truncated metric dimension, we begin by defining $k$-truncated resolving sets.

\begin{definition}
($k$-Truncated Resolving Set.) Let $G=(V,E)$ be a graph. $R \subseteq V$ is a $k$-truncated resolving set of $G$ when, for all $u,v\in V$, if $d_k(u,r)=d_k(v,r)$ for all $r \in R$, $u=v$. 
\end{definition}

\begin{definition}
($k$-Truncated Metric Dimension.) The $k$-truncated metric dimension of a graph $G=(V,E)$ is the size of smallest $k$-truncated resolving sets of $G$ and is denoted $\beta_k(G)$.
\end{definition}

Given a graph $G=(V,E)$ and a $k$-truncated resolving set $R$ on $G$, we define 
\[d_k(v|R) := (d_k(v,r))_{r \in R}.\] 
The \textit{$k$-neighborhood of a vertex} $u$ is the set of all $v\in V$ such that $d(u,v)\le k$. Intuitively, truncated metric dimension restricts the view of each vertex in a graph to its $k$-neighborhood, making all vertices beyond this neighborhood appear identical.

While the results of this paper focus on finite, connected, simple graphs, it is sometimes convenient to focus on matrices underlying the graph structure. 
To that end, the definitions above are easily extended to apply to more general matrices adopting the perspective of multilateration~\cite{tillquist2019low}.

\begin{definition}
($k$-Truncated Resolving Set of a Matrix.) Let $\Mbf$ be an $(n \times m)$ matrix with real entries, and $\Mbf_k$ be the matrix such that $\Mbf_k(a,b) := \min\{\Mbf(a,b), k+1\}$. $R \subseteq \{1,\dots,m\}$ is a $k$-truncated resolving set on $\Mbf$ if, for all distinct $u,v \in \{1,\dots,n\}$, there is an $r \in R$ such that $\Mbf_k(u,r) \neq \Mbf_k(v,r)$.
\end{definition}

\begin{definition}
($k$-Truncated Metric Dimension of a Matrix.) For an $(n \times m)$ matrix $\Mbf$ with real entries, the $k$-truncated metric dimension of $\Mbf$ is $\beta_k(\Mbf) := \beta(\Mbf_k)$. If all rows of $\Mbf$ are equal, its metric dimension is defined as $+\infty$.
\end{definition}

Importantly, for any graph $G$ and any $k \geq 0$, $\beta_k(G) = \beta_k(\Dbf)$ where $\Dbf$ is the distance matrix of $G$. When $k=0$ and $G$ has $n$ vertices, $\Dbf = \mathbf{1} - \mathbb{I}$ where $\mathbf{1}$ is the $(n \times n)$ matrix of all ones and $\mathbb{I}$ is the $(n \times n)$ identity matrix. As a result, $\beta_0(G) = (n-1)$. \textbf{Unless otherwise stated, we assume $\mathbf{k>0}$.}

Observe that if $\delta = \diam(G)$ is the diameter of $G$, then $\beta_{\delta}(G) = \beta(G)$. The following result establishes a relationship between traditional and truncated metric dimension.

\begin{lemma}
\label{lem:D_min_k}
For any graph $G = (V,E)$, $\beta_k(G)$ is a decreasing function of $k$, namely, $\beta_{k+1}(G) \leq \beta_k(G)$ for all $k>0$. In particular, if $\delta = \diam(G)$, $\beta_k(G) = \beta(G)$ for all $k \geq \delta$.
\end{lemma}

\noindent\textit{Proof.} Let $R$ be a $k$-truncated resolving set of $G$ and suppose that $r \in R$ distinguishes $u,v \in V$, i.e. $d_k(u,r) \neq d_k(v,r)$. In particular, since $d_k(\cdot,\cdot)\le(k+1)$, we may assume without loss of generality that $d_k(u,r)<(k+1)$. Then $d_{k+1}(u,r) = d_k(u,r)$, and $d_{k+1}(v,r) = d_k(v,r)$ or $d_{k+1}(v,r) = d_k(v,r)+1$. In either case, $d_{k+1}(u,r) \neq d_{k+1}(v,r)$, which implies that $R$ is a $(k+1)$-truncated resolving set of $G$ and $\beta_{k+1}(G) \leq \beta_k(G)$.

Finally, since $\delta$ is the maximum shortest path distance between two nodes in $G$, $d(u,v) = d_k(u,v)$ for all $u,v \in V$ and $k \geq \delta$. Hence, $\beta(G) = \beta_k(G)$.\hfill$\Box$ \\

\section{Paths and Cycles}
\label{sec:paths_cycles}

Let $P_n$ and $C_n$ denote the path and cycle graphs with $n$ vertices, respectively. In the context of traditional metric dimension, these structures are easily understood. In fact, $\beta(G) = 1$ if and only if $G$ is isomorphic to $P_n$ for some $n>0$~\cite{chartrand2000resolvability}, and $\beta(C_n) = 2$ for all $n > 2$~\cite{chartrand2000cycles}. For $P_n$, either endpoint will serve as a minimal resolving set, while any pair of distinct vertices $u,v \in C_n$ such that $d(u,v) \neq \frac{n}{2}$ resolves the vertices of a cycle. 

When only truncated distance information is available, the situation is somewhat more complex. The next result establishes a simple but non-intuitive formula for the $k$-truncated metric dimension of paths. It extends a characterization of minimal locating dominating sets~\cite{slater1987domination}, that is $1$-truncated resolving sets such that every vertex is adjacent to at least one element of the set.

\begin{theorem}
\label{thm:Pn_min_k}
Let $n,k>0$ and define $m := \lfloor \frac{n}{3k+2} \rfloor$. Then 
\[\beta_k(P_n) = \begin{cases} 
      1, & \text{if } n=1; \\
      2m, & \text{if } n \bmod (3k+2) \in \{0, 1\} \\
      2m+1, & \text{if } n \bmod (3k+2) \in \{2, \dots, k+2\} \\
      2m+2, & \text{if } n \bmod (3k+2) \in \{k+3, \dots, 3k+1\}
   \end{cases}\]
\end{theorem}

\noindent\textit{Proof.} Consider the path $P_n$ with vertices $v_1, \dots, v_n$ and let $R = R' \cup R''$ where 
\[R' = \{v_i | i \bmod (3k+2) \in \{k+1, 2(k+1)\}, i \leq (3k+2)m\}\]
and
\[R''=\begin{cases}
    \{v_n\}, & \text{ if } n=1 \\
    \{\}, & \text{ if } n \bmod (3k+2) \in \{0,1\} \\
    \{v_n\}, & \text{ if } n \bmod (3k+2) \in \{2, \dots, k+2\} \\
    \{v_{(3k+2)m+k+1}, v_{\min(n, (3k+2)m+2(k+1))}\}, & \text{ if } n \bmod (3k+2) \in \{k+3, \dots, 3k+1\}
\end{cases}\]
We claim that $R$ is a minimal $k$-truncated resolving set of $P_n$.

First, we show that $R$ is a $k$-truncated resolving set. If there are two vertices $v_i$ and $v_j$ such that $i<j$ and $d_k(v_i|R) = d_k(v_j|R)$, then, since these vertices are on a path, $d_k(v_i|R)$ and $d_k(v_j|R)$ are either $(k+1, \dots, k+1)$ or $(k+1, \dots, \ell, \dots, k+1)$ where $\ell \leq k$. By construction, the only vertex at truncated distance $(k+1)$ from all elements of $R$ is $v_{(3k+2)m+1}$ when $n>1$ and $n \bmod (3k+2) \in \{1,k+2\}$. On the other hand, if $d_k(v_i|R) = d_k(v_j|R) = (k+1, \dots, \ell, \dots, k+1)$, there is a vertex $v_h \in R$ such that $i < h < j$, $d(v_i, v_h) = d(v_j, v_h) = \ell \leq k$, and $h \bmod (3k+2) \in \{k+1,2(k+1)\}$. If $h \bmod (3k+2) = (k+1)$, then $v_{\min(n,h+k+1)} \in R$. As a result, $d(v_i, v_{\min(n,h+k+1)}) > d(v_j, v_{\min(n,h+k+1)})$ and $d(v_j, v_{\min(n,h+k+1)}) < (k+1)$ so that $d_k(v_i|R) \neq d_k(v_j|R)$. If $h \bmod (3k+2) = 2(k+1)$ instead, $v_{h-(k+1)} \in R$. Then $d(v_j, v_{h-(k+1)}) > d(v_i, v_{h-(k+1)})$ and $d(v_i, v_{h-(k+1)}) < (k+1)$ so that $d_k(v_i|R) \neq d_k(v_j|R)$. Thus, $d_k(v|R)$ is unique for all vertices $v$ in $P_n$ and $R$ is a $k$-truncated resolving set.

Next, we show that $R$ is of minimum size. To begin, note that $R$, and $R'$ more specifically, divides $P_n$ into consecutive subpaths or blocks $B_1, \dots, B_m$ each of size $3k+2$. In particular, letting $w_i = (i-1)(3k+2)$, $v_{w_i+k+1}, v_{w_i+2(k+1)}  \in R'$ resolve the vertices of $B_i$ and, for all vertices $v \in B_i$ and $r \in R' \setminus \{v_{w_i+k+1}, v_{w_i+2(k+1)}\}$, $d(v, r) > k$. When $n \bmod (3k+2) \neq 0$, $B_{m+1}$ is the final subpath of size $q=n-(3k+2)m$.


Now, we attempted to describe a $k$-truncated resolving set $\hat{R}$ of $P_n$ such that $|\hat{R}| < |R|$. Our general approach will be to show that resolving the vertices in $B_i$ depends on vertices in neighboring blocks when $|\hat{R}| < |R|$ and to notice that this leaves vertices unresolved at the ends of the path. We proceed by cases.


\noindent Case 1: If $n=1$, $|R|=1$ and there can be no smaller resolving set.

\noindent Case 2: If $n \bmod (3k+2) \in \{0,1\}$, $n \geq (3k+2)$ and there must be a block $B_i$ containing at most one element $v_j$ of $\hat{R}$. Further, $B_i$ must contain at least one element of $\hat{R}$, otherwise $d_k(v_{w_i+k+1}|\hat{R}) = d_k(v_{w_i+k+2}|\hat{R}) = (k+1, \dots, k+1)$ and these vertices would be indistinguishable.

Note that $j \not \in \{w_i+1, \dots, w_i+k\}$ as this would cause $d_k(v_{w_i+2k+1}|\hat{R}) = d_k(v_{w_i+2k+2}|\hat{R}) = (k+1, \dots, k+1)$. By a symmetric argument, $j \not \in \{w_i+2k+3, \dots, w_i+3k+2\}$. In addition, $j \not \in \{w_i+k+2,\dots,w_i+2k+1\}$ as this would cause $d_k(v_{j-1}|\hat{R}) = d_k(v_{j+1}|\hat{R}) = (k+1, \dots, k+1)$. So $j$ must be $(w_i+k+1)$ or $(w_i+2k+2)$.

If $j = (w_i+k+1)$, the last vertex of the previous block, $v_{w_{i-1}+3k+2} \in B_{i-1}$, must be in $\hat{R}$. If this were not the case or if there were no previous block, $d_k(v_{w_i+k}|\hat{R}) = d_k(v_{w_i+k+2}|\hat{R})$. Subsequently, $v_h$ with $w_{i-1}+k+1 \leq h \leq w_{i-1}+3k+1$ must be in $\hat{R}$ so that $d_k(v_{w_{i-1}+2k+1}|\hat{R}) \neq (k+1, \dots, k+1)$ and $v_{w_{i-1}+2k+1}$ is distinguished from $v_{w_i+2k+2}$. 
Repeating this argument until $B_1$ is reached, we find that at least two vertices in $B_1$ will be indistinguishable unless at least three vertices in $B_1$ are elements of $\hat{R}$. As a result, $|\hat{R} \cap B_1 \cap \dots \cap B_i| \geq |R \cap B_1 \cap \dots \cap B_i|$. That is, up to $B_i$, $\hat{R}$ must be at least as large as $R$.

Similarly, $j = (w_i+k+1)$ forces $v_{w_{i+1}+1} \in \hat{R}$ since $d_k(v_{w_i+2k+2}|\hat{R}) = d_k(v_{w_i+2k+3}|\hat{R}) = (k+1, \dots, k+1)$ otherwise. This means $v_h$ where $w_{i+1}+2 \leq h \leq w_{i+1}+k+2$ must be in $\hat{R}$ so that $d_k(v_{w_i+3k+2}|\hat{R}) \neq d_k(v_{w_{i+1}+2}|\hat{R})$. 
Following this same reasoning until the last block of the path is reached, either $B_m$ will include at least three elements of $\hat{R}$ or, if $n \bmod (3k+2) = 1$, $v_n \in \hat{R}$. Thus $|\hat{R} \cap B_i \cap \dots \cap B_{m+1}| \geq |R \cap B_i \cap \dots \cap B_{m+1}|$

As a result, $|\hat{R}| \geq |R|$. A symmetric argument shows the same result when $j=w_i+2(k+1)$. 

\noindent Case 3: If $n \bmod (3k+2) \in \{2, \dots, k+2\}$, there are two subcases. Either there is a block $B_i$ where $1 \leq i \leq m$ with at most one element of $\hat{R}$ or $B_{m+1}$ includes no elements of $\hat{R}$. In the first subcase, an argument similar to that in Case 2 shows $|\hat{R}| \geq |R|$.

In the second subcase, note that, if $q = k+2$, $d_k(v_{n-1}|\hat{R}) = d_k(v_n|\hat{R}) = (k+1, \dots, k+1)$ regardless of which vertices from $B_m$ are in $\hat{R}$.

On the other hand, if $q < k+2$, there must be a vertex $v_j \in \hat{R}$ with $(w_m+2k+1+q) \leq j \leq (w_m+3k+2)$ so that $d_k(v_{n-1}|\hat{R}) \neq d_k(v_n|\hat{R})$. Consequently, there must be a vertex $v_i \in \hat{R}$ that distinguishes $v_{j-1}$ and $v_{j+1}$. In particular, $j-(k+1)\leq i < j$ and the minimum value of $i$ is $w_m+k+q$.

Applying logic similar to that followed in Case 2 and working back to $B_1$, this means that there must be two vertices $v_{j'}$ and $v_{i'}$ in $\hat{R}$ such that $(2k+1+q) \leq j' \leq (3k+2)$ and $j'-(k+1) \leq i' < j'$. If $j = w_m+2k+3$, the smallest possible value of $i'$ is $k+q$ and $d_k(v_1|\hat{R}) = d_k(v_n|\hat{R}) = (k+1, \dots, k+1)$. If $j > w_m+2k+3$, $i' > k+q$ and $d_k(v_1|\hat{R}) = d_k(v_2|\hat{R}) = (k+1, \dots, k+1)$. In either case, there must be at least one additional vertex in $\hat{R}$ to resolve $P_n$ and $|\hat{R}| \geq |R|$.


\noindent Case 4: If $n \bmod (3k+2) \in \{k+3, \dots, 3k+1\}$, there are two subcases. Either there is a block $B_i$ where $1 \leq i \leq m$ with at most one element of $\hat{R}$ or $B_{m+1}$ includes at most one element of $\hat{R}$. Once again, following the logic used in Case 2, we find that $|R|-1$ vertices are not enough to resolve the path in the first subcase.

In the second subcase, if $q>2k+2$ and $v_j \in \hat{R}$ with $j\geq (3k+2)m+2$, either $d_k(v_{n-1}|\hat{R}) = d_k(v_n|\hat{R}) = (k+1, \dots, k+1)$ or $d_k(v_{j-1}|\hat{R}) = d_k(v_{j+1}|\hat{R}) = (k+1, \dots, k+1)$.

Instead, if $q\leq 2k+2$ and $v_j \in \hat{R}$ with $(3k+2)m+2 \leq j < n$, note that there must be a $v_i \in \hat{R}$ such that $d(v_i, v_j) \leq k$ so that $d_k(v_{j-1}|\hat{R}) \neq d_k(v_{j+1}|\hat{R})$. In particular, $j-(k+1) \leq i < j$. As a result, there must be a vertex $v_h \in \hat{R}$ with $i-(2k+1)-1 \leq h < i$ so that at most one vertex is associated with the $k$-truncated distance vector $(k+1, \dots, k+1)$. However, unless $i-(k+1) \leq h$, $v_{h-1}$ and $v_{h+1}$ must be distinguished by an additional vertex in $\hat{R}$. Once again following logic similar to that used in Case 2 and working back to $B_1$, we find that more than two vertices in $B_1$ are required to distinguish all vertices in $P_n$.

If $j=n$, $i\geq w_m+2k+3$ and the approach in the second part of Case 3 may be used to show that $|\hat{R}| \geq |R|$. 

So there is no case in which $P_n$ has a $k$-truncated resolving set $\hat{R}$ smaller than $R$. As a result, $\beta_k(P_n) = |R|$ and the result follows. \hfill$\Box$ \\

\begin{figure}[h] 
\centering 
\input{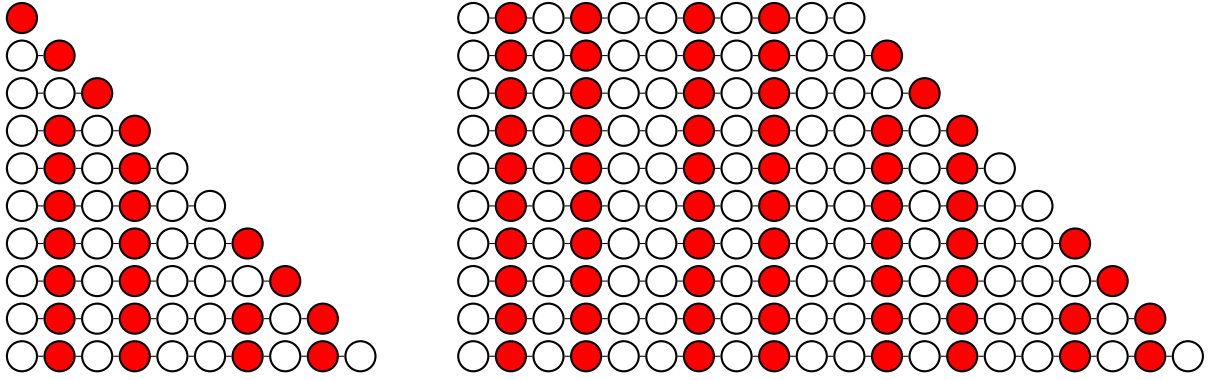}
\caption[$1$-Truncated Resolving Sets on $P_n$]{Minimal $1$-truncated resolving sets of $P_n$ for $1 \leq n \leq 20$ in red.} 
\label{fig:path_min_2_ex_2} 
\end{figure}

Graphs with $k$-truncated metric dimension $1$ are now easily characterized.

\begin{corollary}
\label{cor:Pn_min_k_characterize}
$\beta_k(G) = 1$ if and only if $G \cong P_n$ for some $n \leq k+2$.
\end{corollary}

\noindent\textit{Proof.} Due to Theorem~\ref{thm:Pn_min_k}, $\beta_k(P_n) = 1$ if and only if $n\le k+2$. On the other hand, due to Lemma~\ref{lem:D_min_k}, if $\beta_k(G)=1$ then $\beta(G)=1$, which implies that $G \cong P_n$~\cite[Theorem 2]{chartrand2000resolvability}, from which the corollary follows.\hfill$\Box$ \\\\

Thinking of cycles as paths where the endpoints are adjacent allows for a description of $\beta_k(C_n)$ in terms of $\beta_k(P_n)$. The critical observation is that minimal $k$-truncated resolving sets on different blocks 
of $C_n$ can be made independent 
of one another in the sense that elements of a $k$-truncated resolving set in one block do not affect the number of elements required to resolve any other block.

\begin{corollary}
\label{cor:Cn_min_k}
For $n > 2$ and $k > 0$,
\[\beta_i(C_n) = \begin{cases} 
      2, & n \leq 3k+3 \\
      \beta_k(P_n), & \text{otherwise} 
   \end{cases}
\]
\end{corollary}

\noindent\textit{Proof.} Let $C_n$ be a cycle with $2 < n \leq 3k+3$. Since $\beta(C_n) = 2$, Lemma~\ref{lem:D_min_k} implies that $\beta_k(C_n) \geq 2$. Let $R = \{u,v\}$ such that $C_n$ is divided into two sections $S_1$ and $S_2$ where, without loss of generality, $0 \leq |S_1| \leq k$ and $|S_1|+1 \leq |S_2| \leq 2k+1$. Notice, the size of $S_1$ implies that there is no vertex $w \in S_1$ such that $d(u,w) = k+1$ or $d(v,w) = k+1$. The size of $S_2$ implies that there is at most one vertex $w \in S_2$ such that $d(u,w)=d(v,w)=k+1$. Furthermore, since $|S_1| < |S_2|$, there are no vertices $t,w \in C_n$ such that $d(u,t)=d(u,w)$ and $d(v,t)=d(v,w)$. Thus, $R$ is a resolving set and $\beta_k(C_n) = 2$ in this case.

Next, let $C_n = v_1, \dots, v_n$ be a cycle where $n > 3k+3$ and let $m = \lfloor \frac{n}{3k+2} \rfloor$. Note that $\beta_k(P_{(3k+2)m+1}) = \beta_k(P_{(3k+2)m+1-2k})$. In particular, if $v_{(3k+2)m+2}$ and $v_n$ are elements of a minimal $k$-truncated resolving set of $C_n \setminus \{v_1, \dots, v_{(3k+2)m+1}\}$, at most $2k$ vertices of $\{v_1, \dots, v_{(3k+2)m+1}\}$ are distinguished by this resolving set as well. However, this does not change the number of additional vertices required to resolve $\{v_1, \dots, v_{(3k+2)m+1}\}$.

Furthermore, observe that, for any minimal $k$-truncated resolving set $R'$ of $\{v_1, \dots, v_{(3k+2)m+1}\}$ with $n>(3k+2)m+1$, we have $d(v_{(3k+2)m+2}, r) > k$ and $d(v_n, r) > k$ for all $r \in R'$.
If this were not the case, at least two vertices $u,v \in \{v_1, \dots, v_{(3k+2)m+1}\}$ would have $d_k(v|\hat{R}) = d_k(u|\hat{R})$. As a result,
\[\beta_k(C_n) = \beta_k(P_{(3k+2)m+1}) + \begin{cases}
                                            0, & \text{ if } n-(3k+2)m-1=0 \\
                                            1, & \text{ if } n-(3k+2)m-1 \in \{1, \dots, k+1\} \\
                                            2, & \text{ if } n-(3k+2)m-1 \in \{k+2, \dots, 3k+1\}\\
                                          \end{cases}\]

The second term in the equation above is equivalent to $\beta_k(P_{n-(3k+2)m-1})$ when all vertices are forced to be within distance $(k+1)$ of an element of a $k$-truncated resolving set. In particular, while the first term allows a single vertex $v$ with $d_k(v|\hat{R}) = (k+1, \dots, k+1)$, the second term does not. But then, all together, this is the same as the size of a $k$-truncated resolving set for $P_n$ and the result follows.\hfill$\Box$ \\




\section{Extreme Structures}
\label{sec:extreme}

In this section we show that graph structures for which $k$-truncated metric dimension is maximal are intuitive and align with similar results concerning traditional metric dimension. On the contrary, structures for which $k$-truncated metric dimension is minimal are more complex. We are able to fully characterize graphs of order $n$ with $k$-truncated metric dimension $(n-1)$ and $(n-2)$ and we begin to establish connections between extreme values of $\beta_k$ and graph diameter.


\subsection{Graphs with $k$-Truncated Metric Dimension $(n-1)$ and $(n-2)$}

The complete graph on $n$ vertices has the largest possible metric dimension of all such structures under the traditional definition. Indeed, when $G$ has $n$ vertices, $\beta(G) = (n-1)$ if and only if $G \cong K_n$, i.e. $G$ is the complete graph on $n$ vertices~\cite{chartrand2000resolvability}. The following lemma and corollary show that this is also the case when distances are truncated.

\begin{lemma}
\label{lem:Kn_min_2}
For $G=(V,E)$ connected and not complete with $n = |V| > 2$, $\beta_1(G) < (n-1)$.
\end{lemma}

\noindent\textit{Proof.} Let $G=(V,E)$ be a connected but incomplete graph with $n > 2$ nodes. In particular, there must be $u,v \in V$ such that $\{u,v\} \not \in E$. Furthermore, let $w \in V$ such that $\{v,w\} \in E$. Consider the set $R = V \setminus \{u,w\}$. Since $d_1(v,w)=1$ and $d_1(v,u)=2$, $v \in R$ distinguishes $u$ and $w$. Because all other vertices of $G$ are elements of $R$, $R$ is a $1$-truncated resolving set on $G$ of size $(n-2)$ and $\beta_1(G) < (n-1)$.\hfill$\Box$ \\

\begin{corollary}
\label{cor:Kn_min_k}
For $n \geq 2$ and $k > 0$, $\beta_k(G) = (n-1)$ if and only if $G \cong K_n$.
\end{corollary}

\noindent\textit{Proof.} Suppose $G=(V,E)$ is a graph with $n \geq 2$ nodes and $\beta_k(G) = (n-1)$ for some $k > 0$; in particular, from Lemma~\ref{lem:D_min_k}, $(n-1)\leq \beta_1(G)$. Furthermore, $\beta_0(G) = (n-1)$ for any graph, so $\beta_1(G) \leq (n-1)$. Hence, $\beta_1(G) = (n-1)$. Finally, the contrapositive of Lemma~\ref{lem:Kn_min_2} implies that $G$ is complete.

Next, suppose $G$ is a complete graph. Then $\beta(G) = (n-1)$ and its diameter is $\delta = 1$. Lemma~\ref{lem:D_min_k} implies that $\beta_1(G) = \beta(G) = (n-1)$, as claimed.\hfill$\Box$ \\

We now characterize under what conditions $\beta_k(G) = (n-2)$. The result is nearly identical to the associated characterization under traditional metric dimension~\cite[Theorem 4]{chartrand2000resolvability}. In what follows, $G \cup H$ is the disjoint union of the graphs $G$ and $H$, and $G+H$ is the disjoint union of the graphs with all possible edges between $G$ and $H$ added.

\begin{lemma}
\label{lem:characterize_n_minus_2_min_k}
For $k \geq 1$, $\beta_k(G) = (n-2)$ when $n \geq 4$ if and only if $G$ is $K_{s,t}$ ($s,t \geq 1$), $K_s+\overline{K_t}$ ($s\geq 1$ and $t \geq 2$), $K_s + (K_1 \cup K_t)$ ($s, t \geq 1$), or $P_4$ when $k=1$.
\end{lemma}

\noindent\textit{Proof.} Let $G=(V,E)$ be a graph. If $\diam(G) \leq k+1$ the result follows from Theorem 4 in~\cite{chartrand2000resolvability} and Lemma~\ref{lem:D_min_k}. In particular, since the diameters of $K_{s,t}$ ($s,t \geq 1$), $K_s+\overline{K_t}$ ($s\geq 1$ and $t \geq 2$), and $K_s + (K_1 \cup K_t)$ ($s, t \geq 1$) are all 2, $\beta_k$ with $k \geq 1$ for these graphs is $(n-2)$. In addition, if $G \cong P_4$ and $k=1$, we have $\beta_1(P_4) = 2$ from Theorem~\ref{thm:Pn_min_k}. It is left to show that $\beta_k(G) < (n-2)$ with $k \geq 1$ for all $G$ with $\diam(G)>2$ except $P_4$.

Suppose $\diam(G) > k+1$, $k \geq 1$, and $G$ is not $P_4$. Then there must be $v_1,v_4 \in V$ with $d(v_1,v_4) = 3$ and associated shortest path $P = v_1 v_2 v_3 v_4$. We claim that $R_1 = V \setminus \{v_2,v_3,v_4\}$, $R_2 = V \setminus \{v_1,v_3,v_4\}$, or $R_3 = V \setminus \{v_1,v_2,v_4\}$ must be a $k$-truncated resolving set of $G$. To show that this is the case we need only prove that at least one of $R_1$, $R_2$, and $R_3$ differentiates the vertices of $P$.

Let $v \in V \setminus P$ be adjacent to at least one vertex in $P$. Since $P$ is a shortest path between $v_1$ and $v_4$, $v$ may be adjacent to at most 3 elements of $P$ but may not be adjacent to both $v_1$ and $v_4$. We proceed by cases.

\noindent Case 1: Suppose $\{v,v_1\} \in E$. Then $v, v_2 \in R_2$ distinguish the vertices of $P$. In particular, $v$ and $v_2$ define the following representation of vertices in $P$: $v_1 = (1,1)$, $v_2 = (2,0)$, $v_3 = (d_k(v,v_3), 1)$, $v_4 = (d_k(v,v_4), 2)$. If $\{v,v_4\} \in E$ instead, $R_3$ distinguishes the vertices of $P$ by a symmetric argument.

\noindent Case 2: Suppose $\{v,v_2\} \in E$. Then $v, v_3 \in R_3$ distinguish the vertices of $P$. In particular, $v$ and $v_3$ define the following representation of vertices in $P$: $v_1 = (2, 2)$, $v_2 = (1,1)$, $v_3 = (2, 0)$, $v_4 = (d_k(v,v_4), 1)$. If $\{v,v_3\} \in E$ instead, $R_2$ distinguishes the vertices of $P$ by a symmetric argument. 

\noindent Case 3: Suppose $v$ is adjacent to an endpoint of $P$ and its neighbor. Without loss of generality, assume $\{v,v_1\}, \{v,v_2\} \in E$. Then $v, v_3 \in R_3$ distinguish the vertices of $P$. In particular, $v$ and $v_3$ define the following representation of vertices in $P$: $v_1 = (1,2)$, $v_2 = (1,1)$, $v_3 = (2, 0)$, $v_4 = (d_k(v,v_4), 1)$. If $\{v,v_3\}, \{v,v_4\} \in E$ instead, $R_2$ distinguishes the vertices of $P$ by a symmetric argument.

\noindent Case 4: Suppose $v$ is adjacent to an endpoint of $P$ and the vertex at distance two from the endpoint in $P$. Assume that $\{v,v_1\}, \{v,v_3\} \in E$. Then $v, v_1 \in R_1$ distinguish the vertices of $P$. In particular, $v$ and $v_1$ define the following representation of vertices in $P$: $v_1 = (1,0)$, $v_2 = (2,1)$, $v_3 = (1, 2)$, $v_4 = (2, 3)$. If $\{v,v_2\}, \{v,v_4\} \in E$ instead, $v, v_1 \in R_1$ still distinguish the vertices of $P$. In this case we have: $v_1 = (2,0)$, $v_2 = (1,1)$, $v_3 = (2, 2)$, $v_4 = (1, 3)$.

\noindent Case 5: Suppose $\{v,v_2\}, \{v,v_3\} \in E$. Then $v, v_3 \in R_3$ distinguish the vertices of $P$. In particular, $v$ and $v_3$ define the following representation of vertices in $P$: $v_1 = (2, 2)$, $v_2 = (1,1)$, $v_3 = (1, 0)$, $v_4 = (2, 1)$.

\noindent Case 6: Suppose $v$ is adjacent to an endpoint as well as $v_2$ and $v_3$. Assume that $\{v,v_1\}, \{v,v_2\}, \{v,v_3\} \in E$. Then $v, v_3 \in R_3$ distinguish the vertices of $P$. In particular, $v$ and $v_3$ define the following representation of vertices in $P$: $v_1 = (1,2)$, $v_2 = (1,1)$, $v_3 = (1, 0)$, $v_4 = (2, 1)$. If $\{v,v_2\}, \{v,v_3\}, \{v,v_4\} \in E$ instead, $R_2$ distinguishes the vertices of $P$ by a symmetric argument.

Thus, $R_1$, $R_2$, or $R_3$ must serve as a $k$-truncated resolving set of $G$. Furthermore, $|R_1| = |R_2| = |R_3| = (n-3)$. Hence, we have shown that $\beta_k(G) < (n-2)$ with $k \geq 1$ and for any $G$ with $\diam(G) > 2$ except $P_4$. This proves the lemma.\hfill$\Box$ \\

\subsection{Fixed Diameter}

The relationship between traditional metric dimension and graph diameter has been well studied. There are bounds, for example, on the order of trees, outerplanar~\cite{beaudou2018bounding}, interval, and permutation graphs~\cite{foucaud2017identification} with respect to their metric dimension and diameter. The following corollary generalizes a result which states that, for a graph $G$ of order $n$ and with diameter $\delta = \diam(G)$, $n \leq \delta^{\beta(G)}+\beta(G)$~\cite{khuller1996landmarks}.

\begin{corollary}
\label{cor:D_min_k_diam}
Let $G=(V,E)$ and $\beta_k = \beta_k(G)$. Then $|V| \leq (k+1)^{\beta_k} + \beta_k$.
\end{corollary}

\noindent\textit{Proof.} Let $R$ be a $k$-truncated resolving set of $G=(V,E)$ and define $\beta_k := |R|$. For $v \in V \setminus R$, $d_k(v|R) = (d_k(v,r_1), \dots, d_k(v,r_{\beta_k}))$, the vector of distances from $v$ to each element of $R$, has $\beta_k$ dimensions each taking a value in $\{1, \dots, k+1\}$. There are at most $(k+1)^{\beta_k}$ such vectors. For each $v \in R$, $d_k(v|R)$ has exactly one 0, providing $\beta_k$ additional unique vectors. Since there must be at least one distinct vector representation for every $v \in V$, $|V| \leq (k+1)^{\beta_k}+\beta_k$.\hfill$\Box$ \\

Next, we show a simple upper bound on $\beta_k(G)$ given a fixed diameter.

\begin{lemma}
\label{lem:diam_upper_1}
$\beta_k(G) \leq \beta_k(P_{\delta+1}) + (n-(\delta+1))$ for any graph $G=(V,E)$ where $n=|V|$ and $\delta = \diam(G)$.
\end{lemma}

\noindent\textit{Proof.} Let $G=(V,E)$ be a graph with $n=|V|$ and $\delta = \diam(G)$, let $P_{\delta+1} \subseteq V$ be a shortest path in $G$ of length $\delta$ (i.e. order $\delta+1$), and let $R_P \subseteq P_{\delta+1}$ be a $k$-truncated resolving set of $P_{\delta+1}$. Clearly, $R = (V \setminus P_{\delta+1}) \cup R_P$ is an $k$-truncated resolving set of $G$. Then, since $|R| = (n - (\delta + 1))+\beta_k(P_{\delta+1})$, $\beta_k(G) \leq \beta_k(P_{\delta+1}) + (n-(\delta+1))$ as claimed.\hfill$\Box$ \\

\begin{figure}[h] 
\centering 
\begin{tikzpicture}[every node/.style={font=\scriptsize},scale=1.1]
    \node[non-res] (A) at (0,0) {};
    \node[res] (B) at (1,0) {};
    \node[non-res] (C) at (2,0) {};
    \node[res] (D) at (3,0) {};
    \node[non-res,label=c] (E) at (4,0) {};
    \node[non-res] (F) at (4.5,1) {};
    \node[res] (G) at (5,0.5) {};
    \node[res] (H) at (5,-0.5) {};
    \node[res] (I) at (4.5,-1) {};
    
    \draw (A) -- (B) -- (C) -- (D) -- (E) -- (F) -- (G) -- (H) -- (I) -- (F) -- (H);
    \draw (E) -- (G) -- (I) -- (E) -- (H);
\end{tikzpicture}
\caption[$U_{9,5}$]{A visualization of $U_{9,5}$ with a minimal $1$-truncated resolving set in red.} 
\label{fig:U_n_delta_construction} 
\end{figure}
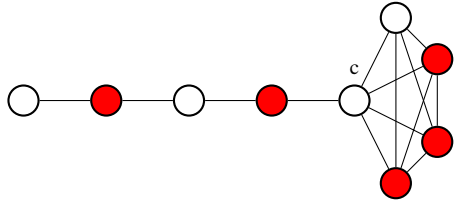

In an effort to tighten this upper bound, we define the following construction. Let $U_{n, \delta}$ with $n \geq \delta$ be a graph consisting of a path $P_{\delta}$ sharing a single endpoint with a complete graph $K_{n-\delta+1}$. In particular, $U_{n, \delta}$ has $n$ vertices and diameter $\delta$. $U_{9,5}$ with a minimal $1$-truncated resolving set is presented as an example in Figure~\ref{fig:U_n_delta_construction}.

\begin{lemma}
\label{lem:U_n_delta_beta_k}
\[ \beta_k(U_{n,\delta}) =
\begin{cases}
    n-\delta, & \text{ if } \delta \leq 2(k+1) \\
    n-\delta+\beta_k(P_{\delta-2k-1}), & \text{otherwise}
\end{cases}.\]
for $n \geq \delta$ and $\delta \geq 1$.
\end{lemma}

\noindent\textit{Proof.} Let $K \subseteq U_{n,\delta}$ and $P \subseteq U_{n,\delta}$ be the complete graph and path subgraphs of $U_{n,\delta}$ respectively and let $\{c\} = K \cap P$. Since the vertices of $K \setminus P$ are indistinguishable with respect to one another and with respect to $P$, any $k$-truncating resolving set $R$ of $U_{n,\delta}$ must include all but one of these vertices. To distinguish $c$ from the unchosen vertex in $K \setminus P$, $R$ must also include $u \in P$ such that $d(u,c) \leq k$. If $\delta \leq 2(k+1)$ and $d(u,c) = \min\{k, \delta-1\}$, $R$ is a minimal $k$-truncating resolving set of $U_{n,\delta}$ of size $(n-\delta)$. 

For $\delta > 2(k+1)$, $d_k(r,v) = k+1$ for all $r \in R$ and for all $v$ at distance at least $2k$ from $c$. Indeed, all vertices within distance $(k+1)$ of $R$ are resolved by $R$ and any remaining unresolved vertices belong to $P$. As a result, the remainder of $P$ may be treated as a path independent from $R$ which has a minimal $k$-truncated resolving set with $\beta_k(P_{\delta-2k-1})$ vertices. The result follows.\hfill$\Box$ \\

It can be shown that $\beta_k(P_{\delta+1})+(n-(\delta+1))-\beta_k(U_{n,\delta}) \in \{0,1\}$. In particular, this construction achieves the bound in Lemma~\ref{lem:diam_upper_1} under certain circumstances.

\section{Trees}
\label{sec:trees}

Many problems that are NP-complete on arbitrary graphs have efficient solutions when restricted to trees. This is true of traditional metric dimension~\cite{harary1976metric,slater1975leaves}. Let $T=(V,E)$ be a tree and call $v \in V$ an \textit{exterior major vertex} if it has degree at least three and there is at least one leaf $u \in V$ for which the path $v,\dots, u$ contains no vertices, except $v$, with degree greater than two. Let $\ell(T)$ be the set of leaves on $T$, $\sigma(T)$ be the set of exterior major vertices on $T$, and $\Delta(v)$ be the set of leaves associated with the exterior major vertex $v \in \sigma(T)$. Then $\beta(T) = |\ell(T)| - |\sigma(T)|$ and $R = \bigcup_{v \in \sigma(T)} \Delta(v) \setminus \{x\}$ is a minimal resolving set where $x$ is any element of $\Delta(v)$~\cite{chartrand2000resolvability}.

As the proof that $R$ is minimal relies on having access to full distance information, the precise relationship between this construction and truncated metric dimension on trees is unclear. However, certain aspects of this proof along with constraints placed on paths by tree structures, suggest that there may be an efficient means of finding minimal $k$-truncated resolving sets on these kinds of graphs. In this section, we present some preliminary results regarding the behavior of truncated metric dimension on trees.

\subsection{The $\mathbb{T}_k$ Family of Trees}

To begin, we define a class of trees for which we can find minimal $k$-truncated resolving sets using the construction described above for traditional metric dimension directly.

Let $\mathbb{T}_k$ be a family of trees defined recursively as follows. Let the empty tree, $T=(\{\}, \{\})$, be in $\mathbb{T}_k$. Then $T=(V,E) \in \mathbb{T}_k$ if four conditions hold.

\begin{enumerate}
    \item $T$ is connected.
    \item There are no vertices of degree two in $V$.
    \item For all minimal non-truncated resolving sets $R$ of $T$ and for all $v \in V$, the vector of distances $d_k(v|R)$ is unique or $(k+1,\dots,k+1)$ (the vector of all $(k+1)$'s).
    \item $T' = T \setminus \{v | \forall u \in V \text{ distinct from } v, d_k(v|R) \neq d_k(u|R)\} \in \mathbb{T}_k$.
\end{enumerate}

Condition $(3)$ may seem difficult to verify at first glance. However, condition $(2)$ significantly restricts the set of minimal non-truncated resolving sets. In particular, condition $(2)$ implies that every minimal non-truncated resolving set $R$ of $T$ must have the form $\bigcup_{v \in \sigma(T)} \Delta(v) \setminus \{x\}$ where $x$ is any element of $\Delta(v)$. 
To see that this is the case, consider an exterior major vertex $v$. Condition $(2)$ implies that $v$ is directly adjacent to each of its associated leaves in $\Delta(v)$. For any $r \in R \setminus \Delta(v)$, that is any element of a resolving set that is not one of these leaves, $d_k(r|\Delta(v))$ is the same. As a result, any resolving set must include all but one of the elements of $\Delta(v)$. Sets of the form $\bigcup_{v \in \sigma(T)} \Delta(v) \setminus \{x\}$ where $x$ is any element of $\Delta(v)$ satisfy this property and are known to be both resolving and minimal~\cite{chartrand2000resolvability}.

The definition of $\mathbb{T}_k$ leads directly to an algorithm for finding minimal $k$-truncated resolving sets for any $T \in \mathbb{T}_k$. Intuitively, we can find minimal $k$-truncated resolving sets of $T \in \mathbb{T}_k$ by constructing a minimal non-truncated resolving set, removing vertices that this set resolves with $k$-truncated distances, and repeating until we are left with a tree with at most one vertex.

\begin{lemma}
\label{lem:good_trees}
Let $T_0 \in \mathbb{T}_k$ and $R_0$ be a minimal non-truncated resolving set on $T_0$. Further, let the sequence of pairs $(T_1, R_1), \dots, (T_n, R_n)$ be generated by repeated application of condition $(4)$ in the definition of $\mathbb{T}_k$. In particular, $T_j = T_{j-1} \setminus \{v | \forall u \in V \text{ distinct from } v, d_k(v|R_{j-1}) \neq d_k(u|R_{j-1})\}$ where $R_{j-1}$ is a minimal non-truncated resolving set of $T_{j-1}$. If $n>0$ and $|T_n| = 1$, $\beta_k(T_0) = \sum_{j=0}^{n-1} |R_j|$. Otherwise, $\beta_k(T_0) = \sum_{j=0}^n |R_j|$.
\end{lemma}

\noindent\textit{Proof.} First, notice that, following an argument similar to the one above, $R_0$ must be of the form $\bigcup_{v \in \sigma(T_0)} \Delta(v) \setminus \{x\}$ where $x$ is any element of $\Delta(v)$.

Next, consider $T_0$ with all vertices in $T_0 \setminus T_j$ resolved. Looking for a minimal non-truncated resolving set of $T_j$, note that vertices from $T_0 \setminus T_j$ provide no benefit, in terms of resolving set size, over picking vertices exclusively from $T_j$. In particular, suppose $\hat{R}_j \subseteq T_0$ is a minimal non-truncated resolving set of $T_j$. Similar to $R_0$, $\hat{R}_j$ must include exactly one vertex from $|\Delta(v)|-1$ of the subtrees rooted at leaves associated with each $v \in \sigma(T_j)$. Let $u \in R_j$ and $t \in \hat{R}_j$ be a descendant of $u$. Since $T_j \in \mathbb{T}_k$, it must be the case that $R_j$ is a $k$-truncated resolving set of every vertex $v \in T_j$ such that $d_k(u,v) \leq k$. However, $\{v| v \in T_j, d_k(t,v) \leq k\} \subseteq \{v| v \in T_j, d_k(u,v)\}$, i.e. the set of vertices in $T_j$ within distance $(k+1)$ of $t$ is a subset of the set of vertices in $T_j$ within distance $(k+1)$ of $u$. Hence, $t$ cannot be used to resolve any vertices not resolved by $u$ with $k$-truncated distances and $R_j$ is minimal. This shows that $\bigcup_{j=0}^{n-1} R_j$ is a minimal $k$-truncated resolving set of $T_{n-1}$.

Finally, there may be a single vertex $v \in T_0$ such that $d_k(v,r) = (k+1)$ for all $r \in \bigcup_{j=0}^{n-1} R_j$. In particular, if $n>0$ and $|T_n| = 1$, $d_k(v,r) = (k+1)$ where $v \in T_n$ for all $r \in \bigcup_{j=0}^{n-1} R_j$. Since, for all other vertices $u \in T_0 \setminus T_n$, there must be at least one $r \in \bigcup_{j=0}^{n-1} R_j$ such that $d_k(u,r) < (k+1)$, this is a unique representation and $\bigcup_{j=0}^{n-1} R_j$ is a $k$-truncated resolving set of $T_0$. Otherwise, $R_n$ is required to differentiate the vertices of $T_n$. The result follows.\hfill$\Box$ \\

In some sense, the definition of $\mathbb{T}_k$ enforces 
independence between subsequent iterations of this process. Since $R_0$ has no effect on $T_j$ for $j>0$, the choice of $R_j$ does not need to involve $R_0$ or any vertex in $T_0 \setminus T_j$. Without this independence property, the interactions between elements of any resolving set can be difficult to characterize. However, by focusing on $k=1$ and limiting the potential influence of each vertex to its immediate neighbors, we can describe a dynamic program guaranteed to find minimal $1$-truncated resolving sets in polynomial time.

\subsection{$\beta_1$ on Trees}

Let $T=(V,E)$ be a tree with an arbitrary root and at least two vertices and let $T_v$ be the subtree of $T$ rooted at $v \in V$. For all $v \in V$ let $C(v)$ be the set of children of $v$. We call $R \subseteq V$ a \textit{locating dominating set} when $R$ is a $1$-truncated resolving set and, for each $v \in V$, there is at least one $r \in R$ such that $d_1(v, r) \leq 1$. Put another way, every $v \in V \setminus R$ must be adjacent to a different non-empty subset of $R$. There is an existing algorithm for finding minimal locating dominating sets on trees, though it is not obvious how this approach might be modified to find minimal $1$-truncated resolving sets~\cite{slater1987domination}. In this section, we describe a novel dynamic programming based algorithm for determining $\beta_1(T)$ exactly on trees.

For all $v \in V$, define the following functions:

\begin{itemize}
    \item $R(v)$ is the size of minimal locating dominating sets for $T_v$. $R_v$ is one such set.
    \item $R'(v)$ is the size of minimal locating dominating sets $R \subseteq T_v$ for $T_v \setminus \{v\}$ such that there is at least one $r \in R$ with $d_k(r,v) \leq k$. $R'_v$ is one such set. 
    \item $R''(v)$ is the size of minimal locating dominating sets $R \subseteq T_v$ for $T_v \setminus \{v\}$. $R''_v$ is one such set.
    \item $R'''(v)$ is the size of minimal locating dominating sets $R \subseteq T_v \setminus \{v\}$ for $T_v$. $R'''_v$ is one such set.
\end{itemize}

It is easy to see that, for all $v \in \ell(T)$ except possibly the root, $R(v) = 1$, $R'(v) = 1$, and $R''(v) = 0$. Note that $R'''(v)$ is undefined for leaves but, as we will see shortly, it can be defined non-recursively. We describe expressions for each of these functions before presenting the algorithm itself.

Assume that we have $R(u)$, $R'(u)$, $R''(u)$, and $R'''(u)$ for all $u \in C(v)$ for some $v \in V$. Consider a locating dominating set $R_v$ of $T_v$. Either $v \in R_v$ or $v \not \in R_v$. In the first case, all children of $v$ are adjacent to at least one element of $R_v$, namely $v$. To guarantee that each child is adjacent to a different non-empty subset of $R_v$, there may be a single $u \in C(v)$ adjacent only to $v \in R_v$ while all other $w \in C(v) \setminus \{u\}$ must be adjacent to at least one other vertex of $R_v$. Consequently, $R_v = \{v\} \cup R''_u \cup (\bigcup_{w \in C(v) \setminus \{u\}} R'_w)$ for some choice of $u \in C(v)$. Then, in this case and taking each $u \in C(v)$ into account, $R(v) = 1 + \min_{u \in C(v)}\{R''(u) + \sum_{w \in C(v) \setminus \{u\}} R'(w)\}$ (equation~\ref{eq:R_v_1}).

If $v \not \in R_v$ instead, there are two additional possibilities. Either at least one or at least two children of $v$ are included in $R_v$. Suppose $u \in C(v)$ is an element of $R_v$. Since it is not guaranteed that $v$ is adjacent to some other vertex in $R_v$, all children of $u$ must be adjacent to at least one other element of $R_v$. The remaining children of $v$ must be located and dominated without $v$. Thus, $R_v = \{u\} \cup (\bigcup_{w \in C(u)} R'_w) \cup (\bigcup_{w \in C(v) \setminus \{u\}} R_w)$ for some choice of $u \in C(v)$ and $R(v) = \min_{u \in C(v)}\{1 + \sum_{w \in C(u)} R'(w) + \sum_{w \in C(v) \setminus \{u\}} R(w)\}$ (equation~\ref{eq:R_v_2}).

Next, suppose $u,w \in C(v)$ are both in $R_v$. Because $v$ is the only vertex that can possibly be adjacent to both $u$ and $w$, we follow an argument identical to when $v \in R_v$ but focus on $u$ and $w$. In particular, to guarantee that each child of $u$ is adjacent to a different non-empty set of $R_v$, there may be a single $x \in C(u)$ adjacent only to $u \in R_v$ while all $y \in C(u) \setminus \{x\}$ must be adjacent to at least one other vertex of $R_v$. A symmetric argument applies to children of $w$. The remaining children of $v$ must be located and dominated without $v$. As a result, $R_v$ can be expressed as the union of the three sets
\begin{align*}
    R_u &= \{u\} \cup R''_{x_u} \cup (\bigcup_{y \in C(u) \setminus \{x_u\}} R'_y) \\
    R_w &= \{w\} \cup R''_{x_w} \cup (\bigcup_{y \in C(w) \setminus \{x_w\}} R'_y) \\
    R_z &= \bigcup_{x \in C(v) \setminus \{u,w\}} R_x
\end{align*}
\noindent for some choice of the pair $u,w \in C(v)$ and for some choice of $x_u$ and $x_w$. Thus, in this case, $R(v)$ can be described with equations~\ref{eq:R_v_3} -- \ref{eq:R_v_5} below and equations~\ref{eq:R_v_1} -- \ref{eq:R_v_5} fully describe $R(v)$.

\begin{align}
    R(v) = \min\{&1 + \min_{u \in C(v)}\{R''(u) + \sum_{w \in C(v) \setminus \{u\}} R'(w)\}, && \label{eq:R_v_1} \\
                 &\min_{u \in C(v)}\{1 + \sum_{w \in C(u)} R'(w) + \sum_{w \in C(v) \setminus \{u\}} R(w)\}, && \label{eq:R_v_2} \\
                 &\min_{u,w \in C(v)}\{2 + \min_{x \in C(u)}\{R''(x) + \sum_{y \in C(u) \setminus \{x\}} R'(y)\} && \label{eq:R_v_3} \\
                 &\phantom{................}+ \min_{x \in C(w)}\{R''(x) + \sum_{y \in C(u) \setminus \{x\}} R'(y)\} && \label{eq:R_v_4} \\
                 &\phantom{................}+ \sum_{x \in C(v) \setminus \{u,w\}}R(x)\}\} && \label{eq:R_v_5}
\end{align}

$R'(v)$ is nearly identical to $R(v)$. The only difference occurs when $v \not \in R'_v$. Since we are not concerned with ensuring that $v$ is adjacent to a different non-empty subset of $R'_v$ as compared to all other vertices, we can focus on the case when at least one child of $v$ is in $R'_v$. Suppose $u \in C(v)$ is in $R'_v$. Following an argument similar to when $v \in R_v$, to guarantee that each child of $u$ is adjacent to a different non-empty subset of $R'_v$, there may be a single $w \in C(u)$ adjacent only to $u$ while all other $x \in C(u)$ must be adjacent to at least one other vertex of $R'_v$. In this case, the remaining children of $v$ must be located and dominated without $v$. This yields equations~\ref{eq:Rp_v_2} and \ref{eq:Rp_v_3} and a full definition of $R'(v)$ below.

\begin{align}
    R'(v) = \min\{&1 + \min_{u \in C(v)}\{R''(u) + \sum_{w \in C(v) \setminus \{u\}} R'(w)\}, && \label{eq:Rp_v_1} \\
                  &\min_{u \in C(v)}\{1 + \min_{w \in C(u)}\{R''(w) + \sum_{x \in C(u) \setminus \{w\}} R'(x)\} && \label{eq:Rp_v_2} \\
                  &\phantom{..............}+ \sum_{w \in C(v) \setminus \{u\}} R(w)\}\} && \label{eq:Rp_v_3}
\end{align}

For $R''(v)$, we do not require that $v$ be adjacent to any element of $R''_v$. However, all children of $v$ must be both located and dominated. So, if $v \not \in R''_v$, we need sets $R_u$ for each $u \in C(v)$ (equation~\ref{eq:Rpp_v_2}). Again, the case when $v \in R''_v$ is identical to the corresponding cases for $R(v)$ and $R'(v)$ (equation~\ref{eq:Rpp_v_1}).

\begin{align}
    R''(v) = \min\{&1 + \min_{u \in C(v)}\{R''(u) + \sum_{w \in C(v) \setminus \{u\}} R'(w)\}, && \label{eq:Rpp_v_1} \\
                   &\sum_{u \in C(v)} R(u)\} && \label{eq:Rpp_v_2}
\end{align}

Finally, $R'''(v)$ follows directly from $R(v)$ when $v \not \in R_v$ (equations~\ref{eq:Rppp_v_1} -- \ref{eq:Rppp_v_4}). We note here that, when $u \in C(v)$, $R'''(u)$ forces $d_1(v,r) = 2$ for all $r \in R'''_u$.

\begin{align}
    R'''(v) = \min\{&\min_{u \in C(v)}\{1 + \sum_{w \in C(u)} R'(w) + \sum_{w \in C(v) \setminus \{u\}} R(w)\}, && \label{eq:Rppp_v_1} \\
                 &\min_{u,w \in C(v)}\{2 + \min_{x \in C(u)}\{R''(x) + \sum_{y \in C(u) \setminus \{x\}} R'(y)\} && \label{eq:Rppp_v_2} \\
                 &\phantom{................}+ \min_{x \in C(w)}\{R''(x) + \sum_{y \in C(u) \setminus \{x\}} R'(y)\} && \label{eq:Rppp_v_3} \\
                 &\phantom{................}+ \sum_{x \in C(v) \setminus \{u,w\}}R(x)\}\} && \label{eq:Rppp_v_4}
\end{align}

We are now ready to define an algorithm for finding minimal $1$-truncated resolving sets on trees.

\begin{algorithm}[!h]
\caption{Minimal $1$-Truncated Resolving Sets on Trees}
\label{beta_2_trees}
\begin{algorithmic}[1]
\Statex{Input: $T=(V,E)$, a tree with $|V|>2$}
\Statex{Output: The minimum size of $1$-truncated resolving sets of $T$}
\Function{$\beta_1$}{$T$}
\State{$S \gets \{R(v)\}$ with any $v \in V$ as the root}
\ForAll{$v \in V$}
\State{$S \gets S \cup \{\sum_{u \in C(v)} R'''(u)\}$}
\EndFor \State{}
\!\!\Return{$\min(S)$}
\EndFunction
\end{algorithmic}
\end{algorithm}

Intuitively, Algorithm~\ref{beta_2_trees} determines the size of minimal locating dominating sets on $T$ and then considers the possibility that each $v \in V$ may be the only vertex not adjacent to any element of a $1$-truncated resolving set. In particular, suppose $v \in V$ is to be this vertex. Since every $u \in C(v)$ must be located and dominated but cannot be included in any $1$-truncated resolving set in this case, we are interested in $R'''_u$ for every $u$. Then, by definition, $R = \bigcup_{u \in C(v)} R'''_u$ locates and dominates every vertex $w \in V \setminus \{v\}$ guaranteeing that $d_1(v,r) = 2$ for all $r \in R$ as desired. So, by taking a minimum of all these values and the size of a minimal locating dominating set on $T$, we are sure to find $\beta_1(T)$. Employing well established methods for keeping track of vertices solving the minimizations in $R(v)$ and $R'''(v)$, Algorithm~\ref{beta_2_trees} can be modified to return a minimal $1$-truncated resolving set of $T$.

\subsection{Extreme Tree Constructions}
\label{sec:trees:extreme}

We end our exploration of truncated metric dimension on trees by examining structures in this family with extreme values of $\beta_k$. As with general graphs in Section~\ref{sec:extreme}, determining a tree structure on $n$ vertices with maximal $k$-truncated metric dimension is straightforward. 

\begin{lemma}
\label{lem:star_tree_max}
For any tree $T$ with $n>2$ vertices, $\beta_k(T) \leq \beta_k(S_n)$ with $k \geq 1$ where $S_n$ is the star graph on $n$ vertices.
\end{lemma}

\noindent\textit{Proof.} To begin, notice that $S_n$ is isomorphic to $K_{1,n-1}$, the complete bipartite graph with partitions of size 1 and $(n-1)$. By Theorem 4 in~\cite{chartrand2000resolvability}, $\beta(S_n) = (n-2)$. Since $K_{1,n-1}$ is not a complete graph for $n>2$, Lemmas~\ref{lem:D_min_k} and~\ref{lem:Kn_min_2} imply $(n-2) \leq \beta_k(S_n) \leq \beta_1(S_n) < (n-1)$ so that $\beta_k(S_n) = (n-2)$ for all $k \geq 1$.

Let $T$ be any tree with $n>2$ vertices. Since $T$ must have $(n-1)$ edges, it cannot be isomorphic to $K_n$. Applying Lemmas~\ref{lem:D_min_k} and~\ref{lem:Kn_min_2} once again, $\beta_k(T) \leq \beta_1(T) < (n-1)$ for all $k\geq 1$ and the result follows.\hfill$\Box$ \\

\begin{figure}[h] 
\centering 
\begin{tikzpicture}[every node/.style={font=\scriptsize},scale=1.1]
    \node[non-res] (A1) at (0,0.5) {};
    \node[non-res] (A2) at (1,0.5) {};
    \node[non-res] (A3) at (2,0.5) {};
    \node[non-res] (A4) at (3,0.5) {};
    \node[res,label=$r_2$] (r2) at (4,0.5) {};
    \node[non-res] (A5) at (5,0.5) {};
    \node[non-res] (A6) at (6,0.5) {};
    \node[non-res] (A7) at (7,0.5) {};
    \node[non-res] (A8) at (5,1.5) {};
    \node[non-res] (A9) at (6,1.5) {};
    \node[non-res] (A10) at (7,1.5) {};
    \node[non-res] (A11) at (6,2.5) {};
    
    \draw (A1) -- (A2) -- (A3) -- (A4) -- (r2) -- (A5) -- (A6) -- (A7);
    \draw (A5) -- (A8);
    \draw (A6) -- (A9) -- (A11);
    \draw (A7) -- (A10);
    
    \node[non-res] (B1) at (0,-0.5) {};
    \node[non-res] (B2) at (1,-0.5) {};
    \node[non-res] (B3) at (2,-0.5) {};
    \node[non-res] (B4) at (3,-0.5) {};
    \node[res,label=$r_3$] (r3) at (4,-0.5) {};
    \node[non-res] (B5) at (5,-0.5) {};
    \node[non-res] (B6) at (6,-0.5) {};
    \node[non-res] (B7) at (7,-0.5) {};
    \node[non-res] (B8) at (5,-1.5) {};
    \node[non-res] (B9) at (6,-1.5) {};
    \node[non-res] (B10) at (7,-1.5) {};
    \node[non-res] (B11) at (6,-2.5) {};
    
    \draw (B1) -- (B2) -- (B3) -- (B4) -- (r3) -- (B5) -- (B6) -- (B7);
    \draw (B5) -- (B8);
    \draw (B6) -- (B9) -- (B11);
    \draw (B7) -- (B10);
    
    \node[res,label=$r_1$] (r1) at (8,0) {};
    \node[non-res] (C1) at (9,0) {};
    \node[non-res] (C2) at (10,0) {};
    \node[non-res] (C3) at (11,0) {};
    \node[non-res] (C4) at (12,0) {};
    \node[non-res] (C5) at (13,0) {};
    
    \draw (A7) -- (r1) -- (B7);
    \draw (r1) -- (C1) -- (C2) -- (C3) -- (C4) -- (C5);
    
\end{tikzpicture}
\caption[$\widetilde{S}_{3, 4}$]{A visualization of $\widetilde{S}_{3, 4}$ with a minimal $4$-truncated resolving set in red.} 
\label{fig:S_tilde_beta_k_construction} 
\end{figure}
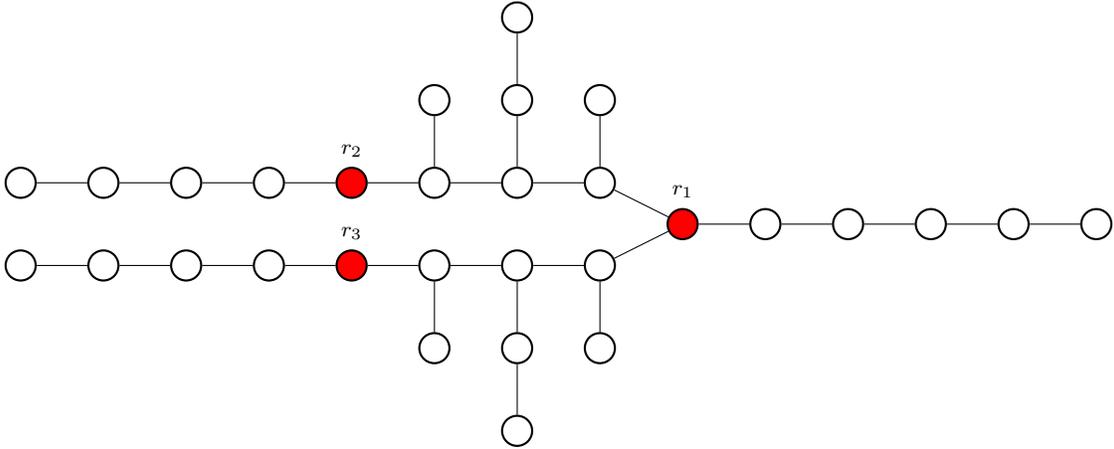

Next, we define a family of trees $\widetilde{S}_{\beta, k}$ such that $\beta_k(\widetilde{S}_{\beta, k}) = \beta$. Let $R = \{r_1, \dots, r_{\beta}\}$ and, for each $r_j \in R$, construct a path of length $k$ with $r_j$ as an endpoint. Include a single extra vertex at the end of the path associated with $r_1$. This vertex will have truncated distance $k+1$ to all elements of $R$. Now, for each $r_j \in R \setminus \{r_1\}$, add a path to $r_1$ of length $(k+1)$ if $k=1$ and length $k$ otherwise with $r_j$ and $r_1$ as endpoints. For each vertex $v$ on a path connecting $r_1$ and $r_j$, add a new path of length $k-\max\{d(v, r_1), d(v, r_j)\}$ with $v$ as an endpoint. In particular, the other endpoint of these paths will be at distance $k$ from at least one of $r_1$ and $r_j$. The number of vertices added for these paths is $\lfloor k^2/4 \rfloor$. As a result, $\widetilde{S}_{\beta, k}$ has 
\[\widetilde{s}(\beta, k) = 1 + \beta (k+1) + (\beta-1) (k-1 + \lfloor \frac{k^2}{4} \rfloor)\]
vertices. $\widetilde{S}_{3, 4}$ is given as an example in Figure~\ref{fig:S_tilde_beta_k_construction}

Observe that $R$ is a $k$-truncated resolving set of $\widetilde{S}_{\beta, k}$. Indeed, since $d(r_i, r_j) = k$ for each distinct pair $r_i, r_j \in R$, each individual element of $R$ resolves its associated path of length $k$ (or $k+1$ for $r_1$) while $r_j \in R \setminus \{r_1\}$ and $r_1$ together resolve all vertices $v$ such that $0 < d(v, r_1), d(v, r_j) \leq k$. 

To see that $\beta_k(\widetilde{S}_{\beta, k}) = \beta$, note that, for any set of vertices $R'$ in $\widetilde{S}_{\beta, k}$ such that $|R'| < \beta$, there must be at least two vertices $u$ and $v$ with $d_k(u|R') = d_k(v|R') = (k+1, \dots, k+1)$. Thus, $\beta_k(\widetilde{S}_{\beta, k}) \geq \beta$. Furthermore, we believe, though do not show, that $\beta_k(\widetilde{S}_{\beta, k}) \leq \beta_k(T)$ for any tree $T$ with $\widetilde{s}(\beta, k)$ vertices.

\section{Conclusion}

Truncated metric dimension restricts the ability of individual vertices to accurately assess distances to far away points in a graph. This variation on the traditional definition forces resolving sets to take a local perspective and has the potential to provide more useful distance constrained resolving sets in a number of real world scenarios. 

In this work, we introduced this concept and explored connections to the traditional definition as well as behavior on paths, cycles, and certain types of trees. We also investigated graph constructions achieving upper and lower bounds in different circumstances.

A variety of interesting questions remain open. For instance, can the $k$-truncated metric dimension of arbitrary trees be determined efficiently? Can approximations of $k$-truncated metric dimension be used to effectively approximate traditional metric dimension through Lemma~\ref{lem:D_min_k}? How effective a tool are $k$-truncated resolving sets for mitigating problems associated with the accumulation of variance in transmission networks in different types of applications? There are many avenues for future exploration related to these ideas.

\section{Acknowledgements}


This research was partially funded by NSF IIS grant 1836914. 

\bibliographystyle{amsrefs}
\bibliography{biblio.bib}

\end{document}